\documentclass[oupthm]{CUP-JNL-BCM}%

\usepackage{graphicx}
\usepackage{multicol,multirow}
\usepackage{amsmath,amssymb,amsfonts}
\usepackage{mathrsfs}
\usepackage{rotating}
\usepackage{appendix}
\usepackage[numbers]{natbib}
\usepackage{multicol}

\theoremstyle{oupplain}
\newtheorem{theorem}{Theorem}[section]
\newtheorem{lemma}[theorem]{Lemma}
\theoremstyle{oupdefinition}
\newtheorem{definition}{Definition}[section]
\theoremstyle{oupremark}

\theoremstyle{oupproof}
\newtheorem{proof}{Proof}
\newtheorem{proposition}[theorem]{Proposition}

\numberwithin{equation}{section}

\begin{document}

\begin{Frontmatter}

\title[To Symbolic Dynamics Through The Thue-Morse Sequence]{To Symbolic Dynamics Through The Thue-Morse Sequence}

\author{Diyath Pannipitiya}

\authormark{D. N. Pannipitiya}

\address{\orgname{Department of Mathematical Sciences, Indiana University-Purdue University--Indianapolis (IUPUI)}, \orgaddress{\street{Indianapolis}, \state{IN}, \postcode{46202.}}\email{dinepann@iu.edu}}



\abstract{  The celebrated Thue-Morse sequence, or the Prouhet-Thue-Morse sequence (A010060 in the OEIS), has a number of interesting properties and is a rich source to many (counter)examples. We introduce two different square-free sequences on three letters with one of them is equivalent (up-to permutations of letters) to Thue's original square-free sequence on three letters. Then we use them to introduce an explicit method to construct infinitely many number of recurrent points in ${\{0,1\}^{\mathbb{N}}}$ whose orbit closures under the shift map is  minimal, uncountable, and for any two distinct such points their orbit closures are disjoint.}

\end{Frontmatter}

\section{Introduction}

The Thue-Morse sequence is one of the extensively studied sequences in both number theory and dynamical systems. It has been used as a prototype to construct interesting examples and counter-examples in areas such as Symbolic Dynamics and Ergodic Theory \cite{bib4}, \cite{bib5}, \cite{bib6}, \cite{bib7}. In his paper \cite{bib14} Thue constructed a \textit{square-free} sequence (see page 2 for the definition) on three symbols by first constructing such a sequence on four symbols (\textit{Theorem 2.2}). In this paper we will follow a completely different approach to obtain the same sequence (upto permutation of letters). Further, we will introduce a method to construct infinitely many square-free sequences on any $n$-letter alphabet ($n\geq 3$). Then we use such sequences to construct recurrent points in $\{0,1\}^{\mathbb{N}}$ whose orbit closures under the shift map is  minimal and uncountable. We call such points \textbf{Morse-type minimal points} (see page 11).\\

There are several equivalent definitions for The Thue-Morse sequence $M$ (see, for example, \cite{bib1}). We state a few of them without proving the equivalence.  \\

\begin{definition}
    (\textbf{Morphism$\backslash$Substitution Definition}). \\  Define the morphism $\mu$ on the alphabet $\{0,1\}$ by $\mu (0)=01$, $\mu (1)=10$. Then $M$ is the unique fixed point of $\mu$ that begins with 0. In other words,
        $M=\lim_{n\rightarrow \infty}\mu ^{n}(0)$.
\end{definition}

\begin{definition}
    (\textbf{Sum of terms Definition}).\\  Let $M_n$ be the $n^{th}$ term of $M$ where $n\in \mathbb{N} \cup \{0\}$. Then\\ ${M_n=[\text{number of}\ 1's\ \text{in}\ (n)_2] mod\ 2.}$ Here $(n)_2:= n$ written in base 2.
\end{definition}

\begin{definition}
    (\textbf{Recursive Definition}).\\  
    Define $M:=(M_n)_{n=0}^{\infty}$ recursively by $M_0=0,\ M_{2n}=M_n,\ M_{2n+1}=1-M_n$.
\end{definition}

\begin{definition}
    (\textbf{Flipping Definition}). \\
    Let $M^N:=(M_n)_{n=0}^N$ be the first $N$ terms of $M$ (first $N$-block$\backslash$ $N$-word). Then $M:=\lim_{n\rightarrow \infty}M^{2^n}$ where $M^0=0$ and $M^{2^{n+1}}:=M^{2^{n}}\overline{M^{2^{n}}}$. Here $\overline{M^n}$ is the $n$-block we get by switching $0$ and $1$ in $M^{{n}}$.
\end{definition}

In this paper we use the \textit{Flipping Definition}. \\

Following the terminology of Thue, Morse and Hedlund (\cite{bib9}), we restate the following definitions.

\begin{definition}
     We say a sequence $\textbf{x}$ is 
     \begin{enumerate}
         \item [$\bullet$] \textbf{square-free} (or \textit{irreducible}) if no word of the form $BB$ appears in $\textbf{x}$. 
         \item[$\bullet$] \textbf{cube-free} if no word of the form $BBB$ appears in $\textbf{x}$.
          \item[$\bullet$] \textbf{overlap-free} if no word of the form $BBb$ appears in $\textbf{x}$, where $B$ is a block (word) of any size in $x$, and $b$ is its first letter.
     \end{enumerate}
     
\end{definition}

Notice that the overlap-free property implies the cube-free property. \\

First let us state the conventions and notations that we will be using throughout this paper. Let $M$ be the Thue–Morse Sequence with \textit{Flipping Definition} as mentioned before. For the convenience, let's use the symbol $``|"$ for the \textit{cuts} (or \textit{mirrors}) that \textit{arise} from the \textit{Flipping Definition}. That is, think of $\overline{M^{2^{n}}}$ in $M^{2^{n+1}}$ as the \textit{mirror image} of $M^{2^n}$ through $``|"$. Then 
\begin{align*}
{\small M=0|1|10|1001|1001,0110|1001,0110,0110,1001|1001,0110,0110,1001,0110,\cdots}
\end{align*}
For a fixed $k\in \mathbb{N}$, we say a block $B$ of length $2^k$ in $M$ is a \textit{partition element} (or a \textit{partition block}) if its starting letter is the $(2^kn)^{th}$ term of $M$ for some $n\in \mathbb{N}\cup \{0\}$. Notice that for any $k\in \mathbb{N}$, there are only two types of partition blocks in $M$; the one which starts with $0$ (let's say $x_k$) and the one which starts with $1$ (let's say $y_k$). It is clear that $x_k$ and $y_k$ are the complements of each other. That is; we can get $y_k$ from $x_k$ (similarly $x_k$ from $y_k$) simply by swapping $0$ and $1$.\\

\textbf{Remarks.}
\begin{enumerate}
\item[$I.$] For each $k\in \mathbb{N}$, we call $x_k$ and $y_k$ the \textit{partition blocks} of length $2^k$ in $M$ such that $x_k$ always starts with 0 and $y_k$ always starts with 1.
\item[$II.$] For each $k\in \mathbb{N}$, $x_k$ and $y_k$ are different at each entry (letter).
\item[$III.$] By the construction (inductively), 
\begin{enumerate}
\item[(a).] before each \textit{mirror}, $M$ has the block either $x_ky_k$ or $y_kx_k$. 
\item[(b).] after each \textit{mirror}, $M$ has the block $y_kx_k$.
\end{enumerate}   
The reason for part $(a)$ is: in $M$ we should get $x_{k+1}$ or $y_{k+1}$ before each \textit{mirror}. And the reason for part $(b)$ is: in $M$ we should get $y_{k+1}$ after each \textit{mirror}. This is an inductive process.
\end{enumerate}

 Thus, for example, we can write $M$ in ``\textit{partition forms}'' as follows. Here  by $M(x_k,y_k)$, we mean $M$ expressed in the partition elements $x_k$ and $y_k$.  
\begin{align*}
M=M(x_0,y_0)&:=x_0|y_0|y_0x_0|y_0x_0x_0y_0|y_0x_0x_0y_0,x_0y_0y_0x_0| \cdots\\
=M(x_1,y_1)&:= x_1|y_1|y_1x_1|y_1x_1x_1y_1|y_1x_1x_1y_1,x_1y_1y_1x_1| \cdots \\
.\\
.\\
.\\
=M(x_k,y_k)&:=x_k|y_k|y_kx_k|y_kx_kx_ky_k|y_kx_kx_ky_k,x_ky_ky_kx_k| \cdots
\end{align*}

It is not hard to see that $x_0=0, y_0=1$, and for $k\in \mathbb{N}_0$, $x_{k+1}=x_ky_k$ and $y_{k+1}=y_kx_k$.\\

\textbf{Note:} There are blocks in $M$ which look like $x_k$ or $y_k$ but are not \textit{partition blocks}. For example the middle block of length four of $y_2y_2$ is $0110$, which looks like $x_2$, even though it is not an $x_2$. Therefore we need a new notation for theses blocks too.

\begin{definition}
    We say two blocks $B:=b_1\cdots b_n$ and $C:=c_1\cdots c_n$ are \textbf{similar} if $b_k=c_k$ for all $k=1,\cdots, n$. We say two blocks $B:=b_1\cdots b_n$ and $C:=c_1\cdots c_n$ are \textbf{equal}, if they are similar and, $b_1$ and $c_1$ have the same entry in $M$. In other words, $B$ and $C$ are equal if they are at the ``same place'' in $M$.
\end{definition}

For example: for any $k\in\mathbb{N}$, any two partition blocks $x_k$ are similar but not equal.\\

\textbf{Clarification of Notations.} \\
$\bullet$ By $x'_k$ we denote the blocks of length $2^k$ which are similar but not equal to any of the partition blocks $x_k$. Similarly by $y'_k$ we denote the blocks of length $2^k$ which are similar but not equal to any of the partition blocks $y_k$.\\

Thus, for example, $x_kx_k = x_{k-1}y'_ky_{k-1}$ and $y_ky_k = y_{k-1}x'_kx_{k-1}$ (for any $k\in \mathbb{N}$).

\section{Background and revisiting Thue}

First studied by  Eugène Prouhet in 1851.
    It was discovered independently by several people including Axel Thue (1906), Marston Morse (1921), Max Euwe (1929), and Arshon (1937).\\
    
    The following two results, which we will give proofs for the completeness of this paper, are due to Thue (\cite{bib1}, \cite{bib2}, \cite{bib9}). The proof of the \textit{Theorem 2.1} is classical. And similar versions of the proof can be found in many texts. For the proof of \textit{Theorem 2.2}, we follow a new approach.

\begin{theorem} [\textbf{Thue}]
    M is overlap-free.
\end{theorem}

\begin{proof}
We prove that $M$ does not have a word of the form $BBb$ where $B$ is a block of any size and $b\in \{0,1\}$ is it's first term.\\

Notice that, by construction, $M$ does not have a word of the form $x_kx_kx_k$ or $y_ky_ky_k$ for any $k\in \mathbb{N}_0$. Suppose a word of the for $BBb$ exists in $M$. We have two main cases with few sub-cases.\\

\underline{\textit{Case 1: (Length of B is odd).}} Consider $M$ as $M(x_1,y_1)$. That is, consider the Morse sequence with partition blocks $x_1=01$ and $y_1=10$. \\

$\bullet$ \textit{Case 1. A: [The first B block starts with a 1 (similarly, a 0) which is the end term of an $x_1$ (similarly, the end term of a $y_1$)].} \\

Then the second $B$ block has to start with 0 (1). But the length of $B$ is odd. So the first two terms of the second $B$ block is actually consisted with an $x_0$ ($y_0$). This implies the second term of the first $B$ block is 1 (0). Because the second and the third terms of the first $B$ block is a partition block in $M(x_1,y_1)$, then it has to be $y_1$ ($x_1$). By continuing in this way we get $b=0$ ($b=1$). Which is a contradiction. \\

$\bullet$ \textit{Case 1. B: The first B block starts with an $x_1$ (similarly, with a $y_1$).} \\

Then, because the length of $B$ is odd, the first term of the second $B$ is 0 (1) which is the end term of a $y_1$ ($x_1$). So the very last term of the first $B$ block is the first term 1 (0) of this $y_1$ ($x_1$). Because the second term of the first $B$ block is 1 (0), the second and the third terms of the second $B$ form a $y_1$ ($x_1$). By continuing in this way, we get $b=1$ ($b=0$). Which is a contradiction. Another contradiction we get is, the first term of the second $B$ block is 1 (0).\\

\underline{\textit{Case 2: (Length of B is even).}} Then we have $|B|= m2^n$ where $m,n\in \mathbb{N}$ and $m$ is odd.\\

$\bullet$ \textit{Case 2. A}: $m=1$.\\

Then $|B|=2^n$ for some $n\in \mathbb{N}$ and consider $M$ as $M(x_n,y_n)$. That is, consider the Morse sequence with partition blocks $x_n=01\cdots$ and $y_n=10\cdots$. It is clear that $B\neq x_n$ (and $B\neq y_n$) because if so, then $b=1 (0)$. Also, by the construction of $M(x_n,y_n)$, it is clear that $BB\notin \{x_nx'_n, x'_nx_n, y_ny'_n, y'_ny_n, x_ny'_n, y_nx'n, x'_ny_n, y'_nx_n\}$. Only subtle part is what happens if $BB=x'_nx'_n$ (or $BB=y'_ny'_n$). In this case, the first partition block which intersects with $BB=x'_nx'_n$ ($BB=y'_ny'_n$) has it's end part similar to the beginning part of a $x_n$ ($y_n$). So it has to be a $y_n$ ($x_n$). By continuing this argument we can show that there is a word of the form $y_ny_ny_n$ ($x_nx_nx_n$) in $M(x_n,y_n)$ which (over)covers $BB$. Which is a contradiction. \\
\vspace{1cm}

$\bullet$ \textit{Case 2. B}: $m\neq 1$. \\

Then let $B=B_1\cdots B_m$ where $B_i$ is a block of length $2^n$ for $i\in \{1,\cdots, m\}$. Let $p$ be the distance between the first term of $B_1$ and the first term of the first partition block (either an $x_n$ or a $y_n$) which intersect $BB$. So $p\in \{0,\cdots, 2^n-1\}$. Then, because the lengths of $x_n$ and $y_n$ are also $2^n$, by shifting the block $BBb$ by $p$ number of units to left, we get a word of the form $B'B'b'$ where $b'$ is the first term of $B'$. Then considering $M$ as $M(x_{n+1}, y_{n+1})$, we get a situation as exactly as in \textit{Case 1}.\\

This completes the proof of \textit{Theorem 1}.
\end{proof}

\begin{theorem} [\textbf{Thue}]
    Let $\nu$ be the sequence whose $n^{th}$ term is the number of 1's in between $n-1^{st}$ and $n^{th}$ zero of $M$. So,
    \[v= 2 1 0 2 0 1 2 1 0 1 2 \cdots.\]
    Then $v$ is square-free.
\end{theorem}

\begin{proof}
    We will give a different proof than Thue's as follows.
    \begin{enumerate}
    \item [$\bullet$] Replace the first digit of each word of the form $00$ and $11$ with $2$ in $M$. So $00$ becomes $20$ and $11$ becomes $21$. Then we get the following sequence.
  \begin{align*}
\vartheta:=0|2|10|1202|1201,0210|1201,0212,0210,1202|1201,0212\cdots =:(\vartheta_n)_{n=0}
\end{align*}
Here $\vartheta_n$ is the $n^{th}$ term if the sequence $\vartheta$.
\item [$\bullet$]  Define $w=(w_n)_{n=0}:=(\vartheta_n+2mod 3)_{n=0}$. That is, replace $0,1,2$ in $\vartheta$ with $2,0,1$ respectively to get $w$.
\item[$\bullet$] First we prove $v=w$ in \textit{Lemma 2.3}. 
\item[$\bullet$] Then we prove $\vartheta$ is square-free in \textit{Theorem 3.1}.  
\item[$\bullet$] That will conclude $w$, and therefore $v$, is square-free.
\end{enumerate}
\end{proof}

\begin{lemma}
    Let $v$ and $w$ as above. Then $v=w$.
\end{lemma}
\begin{proof}
    We show that both sequences $v$ and $\vartheta$ can be obtained from $M$ by following exactly the same procedure (up-to permutation of letters). Because permutations are bijections, we then have the claim.\\
    
    First consider the construction of $v=(v_n)_{n=0}$. Notice that the mid part of the word $x_2x_2$ ($00$ of $0110\ 0110$) and the mid part of the word $y_2y_2$ ($0110$ of $1001\ 1001$) produce the letters 0 and 2 for $v$, respectively. And, the mid part of the word $x_2y_2$ ($010$ of $0110\ 1001$) and the mid part of the word $y_2x_2$ ($010$ of $1001\ 0110$) produce the letter 1 for $v$. Using this fact and the `flipping property' of $M$, we see that $v_{2^n-1}$ term (the terms which are analogues of the `mirror terms' in $M$) is $1$ or $2$ depending on $n\in \mathbb{N}$ is odd or even respectively. And the terms $v_j,\ j=2^n,\cdots 2^{n+1}-1$ can be obtained by switching 2 and 0 (while fixing 1) in $v_j,\ j=0,\cdots 2^{n}-1$.\\

    Now consider the construction of $\vartheta=(\vartheta_n)_{n=0}$. Because of the same properties of the words $x_2x_2$ and $y_2y_2$ and the flipping property of $M$ mentioned above, we see that $\vartheta_{2^n-1}$ term (the terms which are analogues of the `mirror terms' in $M$) is $2$ or $0$ depending on $n\in \mathbb{N}$ is odd or even respectively. And, the terms $\vartheta_j,\ j=2^n,\cdots 2^{n+1}-1$ can be obtained by switching 0 and 1 (while fixing 2) in $\vartheta_j,\ j=0,\cdots 2^{n}-1$.\\

    Thus we can obtain $v$ by replacing the letters $0, 1, 2$ in $\vartheta$ with the letters $2, 0, 1$ respectively. But by the definition, that is the sequence $w$. Hence $v=w$.
\end{proof}

Lets prove some properties of $M$ before proving the square-free property of $v$. Some of the terminology we use here maybe different from the `standard ones' (of, for example, \cite{bib2}, \cite{bib9}).

\begin{definition}
    We say two blocks $B$ and $C$ \textbf{intersect} each other if $B$ starts in $C$, and vise versa.   
\end{definition}

\begin{lemma}
Let $B$ and $C$ be blocks in $M$ with $B\neq C$. If $C$ intersects $B$, then $C$ is not similar to $B$.
\end{lemma}
\begin{proof}
Let $B$ and $C$ be blocks in $M$ such that $B\neq C$ and $C$ intersects $B$. Without loss of generality, let's assume that $C$ is right to $B$. That is, the first term of $C$ comes some number of terms after the first term of $B$, in $M$. Suppose $C$ is similar to $B$. Let $n=|B|=|C|$ be the length of $B$ (and $C$). Let $D$ be the \textit{intersection block} of $B$ and $C$. Let $k=|D|$ be it's length. So $k\in \{1,\cdots n-1\}$. Because $C$ is similar to $B$, blocks which are similar to $D$ occur at the beginning of $B$ and at the end of $C$ (because $D$ is the end $k$-block of $B$ and the starting $k$-block of $C$. And, $B$ and $C$ are similar). This implies the existence of a word of the form $DEDED$. Here $B=C=DED$. Let $D=dD'$ and $F=D'E$. Then $DEDED=dD'EdD'EdD'=dFdFdD'$. But the existence of $dFdFd$ is prohibited by \textit{Theorem 2.1}. This contradiction completes the prove.  
\end{proof}

We can use \textit{Lemma 2.4} to identify the places at where $x'_k$'s and $y'_k$'s occur.\\

\begin{lemma}
For any $k\in \mathbb{N}$, $x'_k$ occurs only at the center of the blocks of the form $y_ky_k$. Similarly, $y'_k$ occurs only at the center of the blocks of the form $x_kx_k$. 
\end{lemma}

\begin{proof}
Let $k\in \mathbb{N}$. It is clear that there is an $x'_k$ (similarly a $y'_k$) at the center of the block $y_ky_k$ ($x_kx_k$) because $y_ky_k=y_{k-1}x_{k-1}y_{k-1}x_{k-1}=y_{k-1}x'_kx_{k-1}$ (similarly, because  $x_kx_k=x_{k-1}y_{k-1}x_{k-1}y_{k-1}=x_{k-1}y'_ky_{k-1}$). This implies the existence. For the uniqueness: by $Lemma\ 2.4$, wherever there is an $x_k$ ($y_k$) term, there is no an $x_k'$ ($y_k'$) which intersects it. So $x'_k$ ($y_k'$) occurs only in a $y_ky_k$ ($x_kx_k$) block.  Notice that, by \textit{Theorem 2.1}, just before and after $y_ky_k$ ($x_kx_k$), we must have an $x_k$ (a $y_k$). This together with the existence part above with $Lemma\ 2.4$ say that, there can not be another $x'_k$ ($y'_k$) term in the block $y_ky_k$ ($x_kx_k$). Hence the claim.

\end{proof}

\begin{theorem}
Let $B$ be a block in $M$. If there is a word of the form $BB$ in $M$, then $B=x_k$ (or $y_k$) for some $k\in \mathbb{N}$.
\end{theorem}

\begin{proof}
Without loss of generality, let's assume that $k\in \mathbb{N}$ is the maximum so that $B$ has an $x_k$. The proof for the case $y_k$ follows similar arguments. \\

Let $B=Px_kQ$. Then $BB=[Px_kQ][P'XQ']$ where $P'$, $Q'$ and $X$ are similar to $P$, $Q$ and $x_k$ respectively. First we show $X=x_k$. Suppose $X\neq x_k$. Then  by $Lemma\ 2.5$ we get $X=x'_k$. This says that there is a word $y_ky_k$ in this $BB$. \\

\underline{\textit{Case 1: ($|Q|=|Q'|>0$).}}  \\

Then the second $y_k$ is inside the second $B$ block (because there is at-least one $y_k$ term after the $x_k$ term in the first $B$ block). Thus there is an $x_{k-1}$ term in the second $B$ block (the end half of this $y_k$) just after $X$. This implies there is an $x_{k-1}$ term just after the $x_k$ term of the first $B$ block (because $x_k$ is a partition block, it has to be an $x_{k-1}$. Not  an $x'_{k-1}$). Then by the construction of $M$, there should be a $y_{k-1}$ term and an $x_{k-1}$ term just before and after this $x_{k-1}$ term respectively. But then this $y_{k-1}$ term is the end half of a partition block in $M(x_k,y_k)$. So there is an $x_k$ term just before the $x_k$ term in the first $B$ block. Similarly, we have another $x_k$ term just after this $x_k$ term in the first $B$ block. So there is a word of the form $x_kx_kx_k$ in $M$, which is a contradiction.\\

\underline{\textit{Case 2: ($|Q|=|Q'|=0$).}} \\

If $|P|=|P'|=0$, then we have the claim as $BB=x_kx_k$. So suppose ${|P|=|P'|>0}$. Then just like in the $Case\ 1$, we can show an existence of a word of the form $x_kx_kx_k$ or the existence of a word of the form $x_kx_ky_ky_k$ both of which are  contradictions. \\

Thus we must have $X=x_k$. Then $QP'$, which is the block between two $x_k$'s, consists with a string of partition blocks in $M(x_k,y_k)$. Again, we got few cases which depend on the number $N$ of these partition blocks in $QP'$.\\

\underline{\textit{Case 2.1:} $N\geq 4$.} \\

Then one of the two $B$ blocks have at-least two of these partition blocks in addition to the $x_k$ term which is already there. By $Theorem\ 2.1$, these two additional partition block can not be an $x_kx_k$. Thus it must be an $x_ky_k$ or a $y_kx_k$. So in each $B$, there is a word of the form $x_ky_kx_k$. But this gives a contradiction to our maximal choice of $k$. Because we could have chosen, for example, $x_{k+1}y_k$ or $x_ky_{k+1}$ (depending on the position of this word in $M$) instead of $x_ky_kx_k$.\\

\underline{\textit{Case 2.2:} $N=3$.} \\

Let $QP'=UVW$ where $U,V,W\in \{x_k,y_k\}$. Then we can find a word in one of the forms $Wx_kU$ or $VWx_k$ or $UVW$ or $x_kUV$ inside one of the two $B$ blocks. Which gives a contradiction just like in \textit{Case 1}.\\

\underline{\textit{Case 2.3:} $N=2$.} \\

Let $QP'=UV$ where $U,V\in \{x_k,y_k\}$. Then we can find a word in one of the forms $Vx_kU$ or $x_kUV$ or $UVx_k$ inside one of the two $B$ blocks. Or, a word of the form $UVx_kUVx_kU$ or $Vx_kUVx_kUV$ in $M$ (which covers $B$). All of these are impossible by above $Case\ 1$ and$\backslash $or by $Theorem\ 2.1$.\\

\underline{\textit{Case 2.4:} $N=1$.} \\

Let $QP'=U$. Then, by $Theorem\ 2.1$, $U=y_k$ (because there is no $x_kx_kx_k$ in $M(x_k,y_k)$). This together with the fact that $k\in \mathbb{N}$ is the maximum number such that $B$ has an $x_k$, imply the existence of a word of the form $y_kx_ky_kx_ky_k$. But this is impossible by $Theorem\ 2.1$.\\

Thus the only choice we have left is $N=0$. This implies 
\[|P|=|P'|=|Q|=|Q'|=0.\] Hence $BB=x_kx_k$ and therefore $B=x_k$. 

\end{proof}

\begin{lemma} Any block $B$ of $M$ occurs infinitely many times in $M$.
\end{lemma}
\begin{proof}
Let $B$ a block in $M$. Let $d$ be the distance - the number of digits - from the end of $B$ to the nearest, and to the right, \textit{mirror} $m_1$, say. Then by the construction of $M$, another block $B'$, say, which is as same as $B$ occurs in between the next two nearest mirrors $m_2$ and $m_3$  with $d$ being the distance between the end of $B'$ and $m_3$.
\end{proof}

\begin{definition}
A subset $A\subset \mathbb{N}$ (or $\mathbb{Z}$) is said to be \textbf{relatively dense}  if there is $k\in \mathbb{N}$ such that $\{n, n+1, \cdots, n+k\}\cap A \neq \emptyset$ for any $n\in \mathbb{N}$.
\end{definition}

\begin{definition}
Let $X$ be a compact metric space and let $f:X\rightarrow X$ be a continuous function. A point $x_0\in X$ is said to be \textbf{almost periodic} if for any neighborhood $U\ni x_0$, the set $\{i\in \mathbb{N}: f^i(x_0)\in U\}$ is relatively dense in $\mathbb{N}$.
\end{definition}

 The following theorem, which we will give a compatible version to the paper, and without a proof, is due to Gottschalk \cite{bib7}. General treatments can be found in \cite{bib8}.

\begin{theorem}[Gottschalk]
Let $x_0\in \{0,1\}^{\mathbb{N}}$. Then $\overline{\mathcal{O}_{\sigma}(x_0)}$ is minimal if and only if $x_0$ is almost periodic.
\end{theorem}

An equivalent definition for the $definition\ 2.3$ when $X=\{0,1\}^{\mathbb{N}}$ and $f=\sigma$ is the following.

\begin{definition}
    A sequence is called \textbf{uniformly recurrent} (also called \textbf{positively recurrent} in some contexts) if for every sub-word $w$ of finite length, there is $N\in \mathbb{N}$ such that every $N$-letter sub-word contains $w$.
\end{definition}

From the proof of \textit{Lemma 2.7} it is clear that $M$ is uniformly recurrent. Then by the above $Theorem\ 2.8$ we see that $\overline{\mathcal{O_\sigma}(M)}$, the orbit closure of $M$ under the shift map, is a minimal set. Even though the minimality implies $\overline{\mathcal{O_\sigma}(M)}$ can not have any proper periodic orbits, we shall state the following proposition with a different proof which does not need $Theorem\ 2.8$.

\begin{proposition}
$\overline{\mathcal{O_\sigma}(M)}$ is uncountable and contains no periodic points.
\end{proposition}

\begin{proof}
First we prove the uncountability. For this we make a sub-set of $\overline{\mathcal{O_\sigma}(M)}$ which has an obvious one-to-one correspondence with $\{0,1\}^{\mathbb{N}}$. \\

For the first step, we have two choices namely $a_0=x_0=0$ or $b_0=y_0=1$. Suppose we picked $a_0$. Then for the second step we can choose block $a_1$, say, of minimal length such that $a_0a_1$ ends with a block $x_1$ or, we can choose block $a_2$, say, of minimal length such that  $a_0a_2$ ends with a block $y_1$. Similarly, if we had picked $b_0$, then for the second step we can choose block $b_1$, say, of minimal length such that $b_0b_1$ ends with a block $x_1$ or, we can choose block $b_2$, say, of minimal length such that  $b_0b_2$ ends with a block $y_1$. Therefore by the end of the second step, we have four possible choices for a finite sequence with four terms. Continue this process. Suppose we have $2^n$ number of possible choices (different blocks of finite sizes) in the $n^{th}$ step. Pick any of them. Say $a_n$. Then for the $(n+1)^{st}$ step, we can choose a block $a_{1,\cdots, 1}$, say, of minimal length such that $a_na_{1,\cdots, 1}$ ends with a block of $x_{n+1}$ or, we can choose a block $a_{1,\cdots, 2}$, say, of minimal length such that $a_na_{1,\cdots, 2}$ ends with a block of $y_{n+1}$. Thus each choices in the $n^{th}$ step gives birth to two choices in the $(n+1)^{st}$ step. Therefore by the end of the $(n+1)^{st}$ step we have $2^{n+1}$ choices for a finite sequence of length $2^{n+1}$. In other words, this is essentially a binary tree.\\
 
 By $Lemma\ 2.7$ we can continue this process and get sequences of the form $(x_0 \cdots x_i)_{i=1}^{\infty}$ where $x_i\in \{a_i,b_i\}$ for all $i\in \mathbb{N}$. Then the set of all the sequences we obtain in this way, 
 $\{(x_0 \cdots x_i)_{i=1}^{\infty}: x_i\in \{a_i,b_i\}\}$, 
 is uncountable and clearly a sub-set of $\overline{\mathcal{O_\sigma}(M)}$.\\
 
 To show that $\overline{\mathcal{O_\sigma}(M)}$ has no periodic points: if it has a periodic point, then we can find a block $B$ of $M$ and an increasing sequence $(n_k)_{k\in \mathbb{N}}$ of natural numbers such that the sequence of blocks $\{B\underset{n_k - \text{times}}\cdots B\}_{k\in \mathbb{N}}$ appears in $M$. But this is a contradiction to $Theorem\ 2.1$ for any $n_k\geq 3$. Thus $\overline{\mathcal{O_\sigma}(M)}$ has no periodic points.\\
 
 This completes the proof of the proposition.
 
\end{proof}

\section[Two Square-free Sequences]{Two Square-free Sequences}
Consider the Thue–Morse sequence
\begin{align*}
M=0|1|10|1001|1001,0110|1001,0110,0110,1001|1001,0110,0110,1001,0110,1001,\cdots
\end{align*}

\underline{\textbf{Sequence I.}}\\
Replace the second digit of each word of the form $00$ and $11$ in $M$ with $2$. So $00$ becomes $02$ and $11$ becomes $12$. Then we get the following sequence.
\begin{align*}
\theta=0|1|20|1021|2021,0120|1021,0120,2120,1021|2021,0120,2120,1021,0120,1021,\cdots
\end{align*}

 \underline{\textbf{Sequence II.}}\\
  As mentioned before, let $\vartheta$ be the sequence we get by replacing the first digit of the each word of the form $00$ and $11$ in $M$ with $2$ in $M$. So $00$ becomes $20$ and $11$ becomes $21$. Hence,
  \begin{align*}
\vartheta=0|2|10|1202|1201,0210|1201,0212,0210,1202|1201,0212,0210,1201,0210,1202,\cdots
\end{align*}

\begin{theorem}
    Both $\theta$ and $\vartheta$ are square-free.
\end{theorem}

\begin{proof}
    We prove $\theta$ is square-free. Similar arguments will prove $\vartheta$ is square-free. Let's assume that there is a block of the form $B_{\theta}B_{\theta}$ in $\theta$. Let $B_1B_2$ be the corresponding block in $M$. Here $B_1$ corresponds with the first $B_{\theta}$ block and $B_2$ corresponds with the second $B_{\theta}$ block. First we show that $B_1=B_2$. Notice that, the only time we get $B_1\neq B_2$ is when 
\[xB_{\theta}B_{\theta}  =  x[2\cdots y][2 \cdots y] \] where $x,y\in \{0,1\}$ with $x\neq y$. Then $B_1 = [x\cdots y]$ and $B_2=[y\cdots y]$. But then notice that \[xB_1B_2 = x[x\cdots y][y \cdots y] = x[xy\cdots y][yy \cdots y]. \] The second term of $B_1$ must be $y$ as we can not get a block of the form $xxx$ by $Theorem\ 2.1$. This says that the second term of $B_{\theta}$ is $y$. Which implies that the second term of $B_2$ is also $y$. But now we have a block of the form $yyy$. A contradiction. Thus $B_1=B_2=B$, say.  \\

Then by $Theorem\ 2.6$, $BB=x_kx_k$ (or $y_ky_k$) for some $k\in \mathbb{N}$. Notice that
 \begin{enumerate}
 \item[$\bullet$] For each $k\in \mathbb{N}$, $x_k$ always starts with 0 and ends with 0 or 1.
 \item[$\bullet$] The block $x_kx_k$ does not appear at the beginning of $M(x_k,y_k)$. This together with $Lemma\ 2.4$ say that there is a $y_k$ just before and after $BB$.
 \end{enumerate}
 
 Now consider $y_kBBy_k=y_kx_kx_ky_k$. If $x_k$ ends with 0, then $y_k$ ends with 1. Thus $B_{\theta}B_{\theta}$ looks like $B_{\theta}B_{\theta}=[0, \cdots, 0][2, \cdots, 0]$. Contradiction as the two blocks are different. If $x_k$ ends with 1, then $y_k$ ends with 0 and therefore $B_{\theta}B_{\theta}$ looks like $B_{\theta}B_{\theta}=[2, \cdots, 1][0, \cdots, 1]$. Again a contradiction. Thus $B_{\theta}B_{\theta}$ does not exist.
 Hence $\theta$ is a square-free sequence.
\end{proof}

\section{Morse-type minimal points}

\begin{definition}
Let $X$ be a topological space (usually a compact metric space) and let $f:X\rightarrow X$ be a function (usually continuous). Then we say $X$ is minimal if $X$ does not contain any proper, closed, non-empty, and forward $f$-invariant subset.
\end{definition}

By the \textit{proposition} we have that $\overline{\mathcal{O}_{\sigma}(M)}$, the orbit closure of the Morse sequence under the shift map, is uncountable and has no periodic points. This set is called \textit{The Morse Minimal Set} \cite{bib3} (in some contexts it mention as a subset of $\{0,1\}^{\mathbb{Z}}$ \cite{bib2}). So, inspired by the \textit{proposition}, we define the following.

\begin{definition}
Let $X$ be a topological space and let $f:X\rightarrow X$ be a continuous function. We say a point $x\in X$ is \textbf{Morse-type minimal point} if $\overline{\mathcal{O}_f(x)}$ is minimal and uncountable.
\end{definition}

\begin{definition}
Let $X$ and $f$ be as in the above $definition\ 4.2$. We say two pints $x_1, x_2\in X$ are \textbf{Morse-type equivalent} ($x_1 \thicksim_M x_2$) if both $x_1$ and $x_2$ are Morse-type minimal points and $x_1\in \overline{\mathcal{O}_f(x_2)}$.
\end{definition}

\textbf{Remark.} Let $\mathcal{M}$ be the set of all Morse-type minimal points in $X$. Then it is not hard to see that $\thicksim_M$ is an equivalence relation on $\mathcal{M}$.\\

In this section we let $X$ to be a shift invariant (usually closed) subset of  $\{0,1\}^{\mathbb{N}}$ and $f=\sigma$.
 Because $M\in \mathcal{M}$, $\mathcal{M}\neq \emptyset$. 
 Let $[\mathcal{M}]:=\mathcal{M}/\thicksim_M$ be the quotient set of $\mathcal{M}$ by $\thicksim_M$. In this section we show that $[\mathcal{M}]$ is uncountable by introducing a method ($Method\ A$) to construct infinitely many different (with respect to $\thicksim_M$) Morse-type minimal points.\\
 
  The second method ($Method\ B$) we present here gives a different way (in addition to $Method\ A$) to construct points in $\{0,1\}^{\mathbb{N}}$ whose orbit closure under the shift map is uncountable and does not contain any periodic points.\\

\underline{\textbf{{\large Method A.}}}\\

Let's consider the square-free sequence $\theta$ above. It is clear that the uncountablity  of $\overline{\mathcal{O}_{\sigma}(M)}$, the orbit closure of the Thue-Morse sequence under the shift map, implies the uncountability of $\overline{\mathcal{O}_{\sigma}(\theta)}$. Now replace the letters $0,1,2$ in $\theta$ with the  letters $1,2,3$ respectively to get the sequence $\beta$. So, $\beta\in \{1,2,3\}^{\mathbb{N}}$ is square-free, and $\overline{\mathcal{O}_{\sigma}(\beta)}$ is uncountable. Now construct the sequence $\alpha$ by replacing each digit in $\beta$ with a block of 1's whose length is that digit. And separate these blocks of 1's with a block of one 0. Thus
\[\alpha = 10110111010110101110110111010111011 \cdots\]

It is not hard to see that $\overline{\mathcal{O}_{\sigma}(\alpha)}$ is uncountable and is minimal (because $M$ is uniformly recurrent implies $\alpha$ is uniformly recurrent). Hence $\alpha$ is a Morse-type minimal point.\\

\textbf{Note:} 
\begin{enumerate}
    \item [$\bullet$]  We can use the square-free sequence $\vartheta$ also.
 \item [$\bullet$] Using this method, we can construct infinitely many number of different points in $\{0,1\}^{\mathbb{N}}$ whose orbit closures are minimal, uncountable, and disjoint.  Basically by first replacing $\{1,2,3\}$ with any three different positive integers, and then turning it into a point in $\{0,1\}^{\mathbb{N}}$ by following the above mentioned method. Thus the cardinality of $[\mathcal{M}]$ is infinite.
\end{enumerate}

\underline{\textbf{{\large Method B.}}}\\ 

Let $\big(a_n \big)_{n\in \mathbb{N}}\in \{2,3,4\}^{\mathbb{N}}$ be a square-free sequence. By $0_{a_n}$, we denote the block of $0$'s with length $a_n$. For $n\in \mathbb{N}$, let $A_n$ be the set of all blocks of length $n$, in $1$'s and $2$'s. That is $A_n=\{1,2\}^n$. So
\begin{align*}
A_1 &= \{1,2\}\\
A_2 &= \{22,21,12,11\}\\
A_3 &= \{222, 221, 212, 211, 122, 112, 121, 111\}\\
&\cdot \\
&\cdot \\
&\cdot
\end{align*}

Construct $\kappa \in \{0,1\}^{\mathbb{N}}$ as follows. Start with $A_2=\{22,21,12,11\}$. It has four blocks each of length two (with two digits). We replace each digit with a block of 1's whose length is that digit and in between theses 1's, we put $0_{a_1}$. And in between two of these newly constructed blocks, we put a single 0. For the convenience let's bold this 0. So $A_2$ gives us the first few terms of $\kappa$:
\begin{align}
    \kappa=110_{a_1}11 \textbf{0} 110_{a_1}1 \textbf{0} 10_{a_1}11 \textbf{0} 10_{a_1}1\cdots
\end{align}

We continue the process creating the next part (which will be added to $(2)$ in right after a single 0) by considering \[A_3 = \{222, 221, 212, 211, 122, 112, 121, 111\}.\] It has eight blocks each of length three. We replace each digit with 1 and in between theses 1's, we put $0_{a_1}$ and $0_{a_2}$ in order. And in between two of these newly constructed blocks, we put a single 0 (which will be in bold for the convenience). Then

\begin{align*}
    \kappa=&(110_{a_1}11 \textbf{0} 110_{a_1}1 \textbf{0} 10_{a_1}11 \textbf{0} 10_{a_1}1)\textbf{0}(110_{a_1}110_{a_2}11\textbf{0}110_{a_1}110_{a_2}1\textbf{0}110_{a_1}10_{a_2}11\textbf{0}\\
    &110_{a_1}10_{a_2}1\textbf{0}10_{a_1}110_{a_2}11\textbf{0}10_{a_1}10_{a_2}11\textbf{0}10_{a_1}110_{a_2}1\textbf{0}10_{a_1}10_{a_2}1)\cdots
\end{align*}

Notice that $1$'s are determined by $A_n$'s and non-single $0$'s are determined by $\big(a_n \big)_{n\in \mathbb{N}}$. Continue in this way to get $\kappa$. We claim $\overline{\mathcal{O}_{\sigma}(\kappa)}$ is uncountable and has no periodic points.\\

\underline{Uncountability:}
Let $\mathcal{A}\subset \{0,1\}^{\mathbb{N}}$ be the collection of all sequences whose 1's are determined by the sequences in $\{1,2\}^{\mathbb{N}}$ and whose 0's (put in between 1's in the same way as in $\kappa$) are determined by $\big(a_n \big)_{n\in \mathbb{N}}$. Because $\{1,2\}^{\mathbb{N}}$ is uncountable, so is $\mathcal{A}$. We claim $\mathcal{A} \subset \overline{\mathcal{O}_{\sigma}(\kappa)}$. Let $b=\big(b_n\big)_{n\in \mathbb{N}}\in \mathcal{A}$ be arbitrary. \\

 Let $\epsilon >0$ be given. Then we can find $N\in \mathbb{N}$ such that ${\frac{1}{2^n}< \epsilon}$ for all $n\geq N$. Let $m\in \mathbb{N}$ be such that $m>N$. By the construction of $\kappa$, there is a block $B$, say, in $\kappa$ which is as same as the first $m$ terms of the sequence $b$. That is; $B=[b_1...b_m]$. Let $\kappa = \big( \kappa_n\big)_{n\in \mathbb{N}}$ and let $t\in \mathbb{N}$ be such that $b_1=\kappa_{t+1}$. Then $|\sigma^t(\kappa)-b|\leq \frac{1}{2^m}<\epsilon$. This says $b\in \overline{\mathcal{O}_{\sigma}(\kappa)}$ and therefore $\mathcal{A} \subset \overline{\mathcal{O}_{\sigma}(\kappa)}$. And because $\mathcal{A}$ is uncountable, $\overline{\mathcal{O}_{\sigma}(\kappa)}$ is also uncountable.\\
 
 \underline{No periodic points:}
 Suppose $\overline{\mathcal{O}_{\sigma}(\kappa)}$ has a periodic point. Then we can find a block $B$, say, such that the sequence of blocks 
 $\Big(\underbrace{B\cdots B }_\text{$n$ \# of blocks}\Big)_{n\in \mathbb{N}}$ appears in $\kappa$. Thus we can choose a large enough $j\in \mathbb{N}$ so that the block $C=\underbrace{B\cdots B }_\text{$j$ number of blocks}$ has at-least one $0_{a_m}$ for some $m\in \mathbb{N}$. This is possible as the maximum length of a block of all 1's is two but the length of the words of the form $[B\cdots B]$ is increasing. Now consider the block $CCCC$ in $\kappa$. This is possible, again because the length of the words of the form $[B\cdots B]$ is increasing.\\
 
 \underline{\textit{Case 1: (There is no $0_{a_1}$ term in $C$).}}\\
 
 Then because each $C$ has a $0_{a_m}$ and there are three $C$'s in a row, then there must be a block in $\big(a_n \big)_{n\in \mathbb{N}}$ of the form $AA$. But this is a contradiction as $\big(a_n \big)_{n\in \mathbb{N}}$ is a square-free sequence. \\
 
 \underline{\textit{Case 2: (There is a $0_{a_1}$ term in $C$).}}\\
 
 Then in the block $CCCC$, there are at-least three $0_{a_1}$ terms. By the construction of $\kappa$, we put $0_{a_1}$'s in between the first and the second digits of each block in $A_n$ for all $n\in \mathbb{N}$. Thus
 there is $j\in \mathbb{N}$, say, such that $A_j$ contains (at-least) two similar blocks. But this is a contradiction as $A_n$ consists with $2^n$ number of different blocks of length $n$. 
  Therefore $\overline{\mathcal{O}_{\sigma}(\kappa)}$ has no periodic points.\\
 
 \textbf{Note:} Using this method also, we can construct infinitely many different points in $\{0,1\}^{\mathbb{N}}$ with uncountable orbit closure containing no periodic points. Basically by replacing $\{2,3,4\}$ with any three different integers greater than 1.\\
 
 The following theorem says that the Morse-type minimal points are `everywhere'.

 \begin{theorem}
     The set $\mathcal{M}$ of all Morse-type minimal points is dense in $\{0,1\}^{\mathbb{N}}$.
 \end{theorem}

 \begin{proof}
     Let $a=a_0a_1a_2\cdots a_n\cdots \in \{0,1\}^{\mathbb{N}}$ be such that $\overline{\mathcal{O}_{\sigma}(a)}=\{0,1\}^{\mathbb{N}}$. For ${N,k\in \mathbb{N}}\cup\{0\}$, define $A_N(k):=a_{k}\cdots a_{k+N}$ and $B(k):=\bigcup_{N=0}^\infty \{A_N(k)\}$. Notice that $A_N(k)$ is the first $(N+1)^{st}$ block of $\sigma^{k}(a)$. Because $B(k)$ is countable, $B:=\bigcup_{k=1}^\infty B(k)$ is also countable. Let $B=\{b_m:m\in \mathbb{N}\}$ be an enumeration. For $m\in \mathbb{N}$, define $x_m:=b_m\cdot \overline{b_m} \cdot \overline{b_m}\cdot b_m\cdots\in \{0,1\}^{\mathbb{N}}$, where $\overline{b_m}$ is the complement (switching $0$ and $1$) of $b_m$, and they follow the pattern (the \textit{Flipping Definition}) of the Thue-Morse sequence to form $x_m$. The uniform recurrent property of $M$ implies the uniform recurrent property of $x_m$, for any $m\in \mathbb{N}$. Hence $x_m$ is a Morse-type minimal point for any $m\in \mathbb{N}$. Thus the collection $\bigcup_{m=1}^\infty \{x_m\}$ is a subset of $\mathcal{M}$. Because $\overline{\mathcal{O}_{\sigma}(a)}=\{0,1\}^{\mathbb{N}}$, we have that $\overline{\bigcup_{m=1}^\infty \{x_m\}}=\{0,1\}^{\mathbb{N}}$. This implies $\overline{\mathcal{M}}=\{0,1\}^{\mathbb{N}}$. 
 \end{proof}

\textbf{Remark.} We can extend $\mathcal{M}$ as a subset from $\{0,1\}^{\mathbb{N}}$ to $\{0,1\}^{\mathbb{Z}}$. For example, \[M^*:=\cdots m_n\cdots m_2m_1m_0\cdot m_om_1m_2\cdots m_n \cdots \in \{0,1\}^{\mathbb{Z}},\] where $(m_n)_{n=1}=M$, is a Morse-type minimal point as an element of $\{0,1\}^{\mathbb{Z}}$ (in fact sometimes $\overline{\mathcal{O}_{\sigma}(M^*)}$ is called the Morse minimal set \cite{bib2}).\\

\textbf{Open Question} Can we get a point $x^*\in \{0,1\}^{\mathbb{N}}$ from Method B such that $\overline{\mathcal{O}_{\sigma}(x^*)}$ is not minimal?

\section{Acknowledgements} I would like to thank my advisors Prof. Bruce Kitchens and Prof. Roland Roeder of IUPUI for their guidance, helpful comments, and most importantly for asking a lot of interesting questions which undoubtedly have strengthen this paper. If it were not for them, this paper might not exist.

\begin{Backmatter}

\printaddress

\end{Backmatter}

\end{document}